\documentclass{article}
\usepackage{amsmath}
\usepackage{amsfonts}
\usepackage{amsthm}
\usepackage{amscd}
\usepackage{amssymb}
\usepackage{color}
\usepackage{graphicx}
\usepackage{tikz}

\usepackage[colorlinks=true]{hyperref}
\hypersetup{urlcolor=blue, citecolor=red, linkcolor=blue}

\newtheorem{theorem}{Theorem}[section]
\newtheorem{lemma}[theorem]{Lemma}

\newtheorem{observation}[theorem]{Observation}

\theoremstyle{definition}

\theoremstyle{remark}

\theoremstyle{comment}

\newcommand{\1}{\mathbf{1}}

\def\1{\mathbf{1}}

\def\mju{\mathcal{U}}

\def\f2{\mathbb{F}_2}

\newcommand{\disp}{\displaystyle}

\newcommand\remove[1]{}

\newcommand{\lp}[1]{\left( #1 \right)}

\newcommand{\ls}[1]{\left[ #1 \right]}
\newcommand{\lc}[1]{\left\{ #1 \right\}}
\newcommand{\av}[1]{\left| #1 \right|}
\newcommand{\nm}[1]{\left\| #1 \right\|}

\begin{document}

\title{Images of nowhere differentiable Lipschitz maps of $[0,1]$ into $L_1[0,1]$}

\author{Florin Catrina and Mikhail I.~Ostrovskii\\
Department of Mathematics and Computer Science\\
St. John's University\\
8000 Utopia Parkway\\
Queens, NY 11439, USA\\
E-mails: {\tt catrinaf@stjohns.edu} and {\tt
ostrovsm@stjohns.edu}}

\date{\today}
\maketitle

\begin{large}

\noindent{\bf Abstract.} The main result: for every
sequence $\{\omega_m\}_{m=1}^\infty$ of positive numbers
($\omega_m>0)$ there exists an isometric embedding $F:[0,1]\to
L_1[0,1]$ which is nowhere differentiable, but for each $t\in
[0,1]$ the image $F_t$ is infinitely differentiable on $[0,1]$
with bounds $\max_{x\in[0,1]}|F_t^{(m)}(x)|\le\omega_m$ and has an
analytic extension to the complex plane which is an entire
function.
\medskip

\noindent{\bf Keywords:} Differentiability of a Banach space
valued map, isometric embedding, Lipschitz map
\medskip

\noindent{\bf 2010 Mathematics Subject Classification.} Primary:
46G05; Secondary: 46B22.

\section{Introduction}

It is well-known that there exist nowhere differentiable Lipschitz
maps $F:[0,1]\to L_1[0,1]$ (in this paper we denote the function
in $L_1$ corresponding to $t\in[0,1]$ by $F_t$). Apparently the
first and the simplest example of such a map was constructed by
Clarkson \cite[p.~405]{Cla36} (earlier an example with the same
property in  $\ell_\infty[0,1]$ was constructed by Bochner
\cite{Boc33}).
Clarkson's example is $F_t=\mathbf{1}_{[0,t]}$,
where $\mathbf{1}_A$ is the {\it indicator function} of $A$, that
is:

\[\mathbf{1}_A(x)=\begin{cases} 1 &\hbox{ if }x\in A\\
0 &\hbox{ otherwise}.\end{cases}\] It is worth mentioning that the
map constructed by Clarkson is not only a Lipschitz map, but is an
isometric embedding.

So the analogue of the Lebesgue Theorem (see \cite[p.~11]{RS90})
fails for Lipschitz maps from $[0,1]$ to $L_1[0,1]$. Later the
study initiated by Bochner and Clarkson developed into a direction
in Banach space theory devoted to study Banach spaces for which
such phenomena occur. Recall that one of the equivalent
definitions of a Banach spaces with the Radon-Nikod\'ym property
is: A Banach space $X$ has {\it the Radon-Nikod\'ym property} if
and only if every Lipschitz map $F:\mathbb{R}\to X$ is
differentiable almost everywhere. The corresponding part of the
Banach space theory is reflected in many monographs, see
\cite[Chapter 5]{BL00}, \cite{Bou79}, \cite{Bou83}, \cite{DU77},
\cite[Chapter 2]{Pis16}. Radon-Nikod\'ym property plays an
important role in the theory of metric embeddings, see \cite{CK06,
CK09, LN06}.

This paper is motivated by a question of Alessio Figalli
\cite{Fig06} who asked if one may hope for a differentiability
result for Lipschitz maps $F:[0,1]\to L_1[0,1]$ after adding some
smoothness hypothesis on the functions $F_t=F_t(x)$ as functions
of $x$ (for each fixed $t$). His main interest was the case where
we suppose that $F_t(x)$ is $1$-Lipschitz in $x$ for every $t$.
Unfortunately the answer to this question is negative even under
substantially stronger assumptions. The corresponding example is
the main result of this paper. Namely, we prove the following
theorem.

\begin{theorem}
\label{T:FinSm}
Let $\{\omega_m\}_{m=1}^\infty$ be a sequence of positive numbers
 $(\omega_m>0)$. Then there exists an isometric embedding
$F:[0,1]\to L_1[0,1]$ such that

\begin{itemize}

\item[{\bf (A)}] For each $t\in [0,1]$ the function $F_t(x)$ is
infinitely differentiable on $[0,1]$ with bounds $\max_{x\in[0,1]}
\av{F_t^{(m)}(x)} \le \omega_m$ and has an analytic extension to
the complex plane which is an entire function.

\item[{\bf (B)}] The map $t\mapsto F_t$ does not have points of
differentiability in $[0,1]$.

\end{itemize}
\end{theorem}

\section{Proof of Theorem \ref{T:FinSm}}

The idea of the proof can be described in the following way: we
are going to join the function $F_0={\bf 0}$ (identically $0$ on
$[0,1]$) with the function $F_1={\bf 1}$  (identically $1$ on
$[0,1]$) with a $1$-Lipschitz curve in $L_1[0,1]$ consisting of
functions satisfying the condition {\bf (A)} in such a way that
this curve ``$c$-changes direction'' at binary points of
infinitely many ``generations'' with the same $c>0$. The map
$F:[0,1]\to L_1[0,1]$ which we construct is not only
$1$-Lipschitz, but is also an isometric embedding.

By a {\it binary point} in $[0,1]$ we mean a point representable
by a fraction $\frac{i}{2^n}$, where $n$ is a nonnegative integer
and $i$ is an integer satisfying $0\le i\le 2^n$. The {\it
canonical form} of a positive binary point is $\frac{i}{2^n}$,
where $i$ is odd. The {\it generation} of a binary point is the
number $n$ in its canonical form, and we define
that  $0$ is of generation $0$.

We say that $F$  {\it ``$c$-changes direction''} at a binary point
$i/2^n$ of generation $n\ge 1$ if
\begin{equation}
\label{E:cChDir}
||F_{(i-1)/2^n}+F_{(i+1)/2^{n}}-2F_{i/2^n}||_{L_1}\ge c2^{-n}.
\end{equation}

Using a version of a well-known argument we show that if the
``$c$-change of direction'' condition is satisfied at infinitely
many generations of binary points, then the map does not have
points of differentiability, see Lemma \ref{L:cCD_implND}.

The main content of the construction is to show that the
``$c$-change of direction'' condition at sufficiently rare
generations of binary points is compatible with $F_t$ being
infinitely continuously differentiable functions with suitable
bounds on the restrictions of derivatives to $[0,1]$, and having
analytic extensions to the complex plane. First we describe how to
find the functions $F_t$ for all binary $t$.
 These functions  will be chosen to be polynomials.
  Their construction is done generation by
generation. It is easy to see from the construction below that the
map $t\mapsto F_t$ is monotone in the following sense:
\begin{equation}\label{E:Monotone}
0 \leq s \leq t \leq 1 \ \ \mbox{ implies } \ \ 0 \leq F_s(x) \leq
F_t(x) \leq 1 \ \ \mbox{ for all } \ \  x\in [0, 1].
\end{equation}
For this reason the notion of the ``region bounded by the graphs
of $F_{\frac{i-1}{2^{n}}}$ and $F_{\frac{i+1}{2^{n}}}$'' which we
use below is well-defined.
\medskip

The property that $t\mapsto F_t$  is an isometric embedding is
achieved by choosing the image of the binary point $\disp
\frac{i}{2^{n}}$ ($1 \leq i \leq 2^n-1$ and $i$ is odd) of
generation $n$ to be mapped onto a function $F_{\frac{i}{2^{n}}}$
whose graph divides the region bounded by the graphs of
$F_{\frac{i-1}{2^{n}}}$ and $F_{\frac{i+1}{2^{n}}}$ (observe that
${\frac{i-1}{2^{n}}}$ and ${\frac{i+1}{2^{n}}}$ belong to previous
generations) into two parts of equal area. It is clear that the
continuous extension of this map to non-binary points will be an
isometric embedding of $[0,1]$ into $L_1[0,1]$.

The condition that functions $F_t$ satisfy {\bf (A)} is achieved
in the following way: For several (usually very many) consecutive
generations of binary points we define $F_{\frac{i}{2^{n}}}$
(where $\frac{i}{2^{n}}$ is of generation $n$) by
\begin{equation}\label{E:MidPtDef}
F_{\frac{i}{2^{n}}}=\frac12\left(F_{\frac{i-1}{2^{n}}}+F_{\frac{i+1}{2^{n}}}\right).\end{equation}
Note that since $(i-1)$ and $(i+1)$ are even, the functions on the
right hand side are defined earlier. We show that if we use the
definition \eqref{E:MidPtDef} for sufficiently many generations of
binary points, then for the immediately following generation  we
can define $F_{\frac{i}{2^{n+1}}}$ in such a way that it satisfies
both the ``$c$-change of direction" condition and the conditions
which lead to bounds for the derivatives and the existence of an
analytic extension.

\subsection{Plan for estimating the derivatives and showing the
analyticity}\label{S:Plan}

The bounds on the derivatives of $F_t$ and
analyticity are obtained as follows. For each point $t\in[0,1]$ we
find a sequence $\{p(r,t)\}_{r=1}^\infty$ of binary points
converging to $t$ such that the following conditions hold:
\medskip

\noindent{\bf 1.} The series
\begin{equation}\label{E:Series}F_{p(1,t)}+\sum_{r=1}^\infty
\left(F_{p(r+1,t)}-F_{p(r,t)}\right)\end{equation} converges
uniformly  on each bounded disc in the complex plane (this
statement has a natural meaning - recall that
$F_{\frac{i}{2^{n}}}$ are polynomials). It is clear that
restrictions of the terms of this series to $[0,1]$ form a series
which converges to $F_t$, and so this condition together with the
fact that $F_{\frac{i}{2^{n}}}$ are polynomials, implies that
$F_t$ has an analytic extension to $\mathbb{C}$.\medskip

\noindent{\bf 2.} The derivatives of $\{F_{p(r,t)}\}_{r=1}^\infty$
satisfy

\begin{equation}\label{E:Der1}
|F^{(m)}_{p(1,t)}(x)|=0\quad\forall x\in[0,1],
\end{equation}

\begin{equation}\label{E:DerGen}
|F^{(m)}_{p(r+1,t)}(x)-F^{(m)}_{p(r,t)}(x)|\le
\frac{\omega_m}{2^{r}}\quad\forall x\in[0,1].
\end{equation}

As is well-known, \eqref{E:Der1} and \eqref{E:DerGen} together
with the uniform convergence of \eqref{E:Series}  imply that
$|F^{(m)}_t(x)|\le\omega_m$ for all $m\in\mathbb{N}$ and
$x\in[0,1]$. (In more detail, we repeatedly use the statement: if
the series $\sum u'$ of derivatives of a uniformly convergent
series $\sum u$ is uniformly convergent, then $\sum u'$ converges
to the derivative of the  $(\sum u)'$.)
\medskip

A similar argument in a different context was used by Bernstein in
Constructive Function Theory, see \cite[Theorem 13.20 on
p.~119]{Zyg68}.

\subsection{Construction of an example}

We start by letting $F_0={\bf 0}$ (identically $0$ on $[0,1]$) and
$F_1={\bf 1}$  (identically $1$ on $[0,1]$). After that we define
polynomials (with real coefficients) $F_t$ for all other binary points
$t$ in $[0,1]$ generation by generation.
As will be easy to see from our construction, the map
$t\mapsto F_t$ is monotone in the sense of \eqref{E:Monotone}.

We define $F_{\frac{i}{2^n}}$ at a point
$\frac{i}{2^n}$ of generation $n$ so that its graph bisects the
region between the graphs of $F_{\frac{i-1}{2^n}}$ and
$F_{\frac{i+1}{2^n}}$ into two regions of equal area. This will
ensure the isometry condition.

Specifically,  for every $\frac{i}{2^n}$ of generation $n$ we
employ a polynomial $\eta_\frac{i}{2^n}: [0,1]\to
[0, 1]$ to define
\begin{equation}\label{E:MidPtU}
F_{\frac{i}{2^{n}}}=\eta_\frac{i}{2^n} F_{\frac{i-1}{2^{n}}}+
\lp{1- \eta_\frac{i}{2^n}}F_{\frac{i+1}{2^{n}}}
\end{equation}
so that the graph of the restriction of $F_{\frac{i}{2^n}}$ to
$[0,1]$ splits the region betweens the graphs of the restrictions
of $\disp F_{\frac{i-1}{2^{n}}}$ and $\disp F_{\frac{i+1}{2^{n}}}$
to $[0,1]$ into two regions of equal area. For many generations of
binaries, these polynomials are constant: $\disp
\eta_\frac{i}{2^n}(x) = \frac 12$. In this case \eqref{E:MidPtU}
is just \eqref{E:MidPtDef}. However, with this choice of $\eta$
there is no change of direction at $\frac{i}{2^n}$ since in this
case
\[ ||F_{(i-1)/2^n}+F_{(i+1)/2^{n}}-2F_{i/2^n}||_{L_1} = 0.\]

To introduce the  ``$c$-change of direction'' (with $c = 1$) at an
infinite sequence of generations $\{n_r\}_{r=1}^\infty$ we
construct functions $\eta$ with the help of the following Lemma.

\begin{lemma}
\label{L:Split}
 Let $G$ be a polynomial of degree $d\geq 0$
with  $G(x) \geq 0$ for all $x\in [0,1]$. Then there exists
 a polynomial  $\mju$ of degree $(d+1)$  with
 $\mju(x) \in [0, 1]$ for all $x\in [0,1]$  such that the degree $(2d+1)$
polynomials $ \hat G := \mju G$ and $\tilde G:=(1-\mju)G$ satisfy
 \[   ||\hat G||_{L_1} = ||\tilde G||_{L_1} =  \frac12 || G||_{L_1}
 \ \mbox{  and } \
||\hat G-\tilde G||_{L_1} =\frac12  ||G||_{L_1}.\]

\end{lemma}
\proof
Assume first that $||G||_{L_1} = 1$, i.e. $G$ is the density of a probability
distribution on $[0, 1]$. Then the  cumulative distribution
function $\mju(x):=\int_0^x G(t) \ dt$ is a polynomial of degree
$(d+1)$  that maps $[0,1]$ onto $[0,1]$.
Since both $\mju$ and $G$ are nonnegative on $[0,1]$ we
have
\[  ||\hat G||_{L_1} = \int_0^1 \mju (x) G(x) \ dx = \frac12 \int_0^1 \frac{d}{dx}\, \mju^2(x) \ dx
=  \frac12 \,\mju^2(1)-\frac12 \,\mju^2(0) = \frac12.\] Similarly,
\[  ||\tilde G||_{L_1} = \int_0^1 (1-\mju (x)) G(x) \ dx =  \int_0^1 G(x) \ dx -  \int_0^1 \mju (x) G(x) \ dx =\frac12.\]
Note that since $G$ is a polynomial, there exists a unique $m \in
(0,1)$ such that
\[ \mju(m) = \int_0^m G(x) \ dx = \int_m^1 G(x) \ dx =\frac12.\]
Then
\[ \begin{aligned}
& ||\hat G-\tilde G||_{L_1}  = \int_0^1 |1-2\mju (x)| G(x) \ dx = \int_0^m (1-2\mju (x)) G(x) \ dx \\
& + \int_m^1 (-1+2\mju (x)) G(x) \ dx = \int_0^m G(x) \ dx -\int_0^m \frac{d}{dx}\, \mju^2(x) \ dx \\
& -\int_m^1 G(x) \ dx +\int_m^1 \frac{d}{dx}\, \mju^2(x) \ dx =
\frac12 -  \frac14-\frac12 +1-\frac14 = \frac12.
\end{aligned}\]
In the case where  $||G||_{L_1}> 0$ is not necessarily
one, we apply the previous construction to the probability density
$\displaystyle \frac{G}{||G||_{L_1} }$ on $[0,1]$.
\endproof

At certain, rather rare sequence $\{n_r\}_{r=1}^\infty$ of
generations, instead of defining $F_{\frac{i}{2^n}}$ using
\eqref{E:MidPtDef} we do the following. Let $i$ be odd, $1\le i\le
2^{n_r}-1$. Define
\[ G := F_{\frac{i+1}{2^{n_r}}} -F_{\frac{i-1}{2^{n_r}}}.\]
Then $||G||_{L_1} = \frac{2}{2^{n_r}}$. Since
$\frac{i-1}{2^{n_r}}$ and $\frac{i+1}{2^{n_r}}$ are binary of
generations strictly smaller than $n_r$, $F_{\frac{i+1}{2^{n_r}}}$
and $F_{\frac{i-1}{2^{n_r}}}$ are already defined polynomials, and
the monotonicity condition \eqref{E:Monotone} guarantees that $G$
is nonnegative on $[0,1]$. We apply Lemma~\ref{L:Split} to $G$ to
obtain $\mju$ and we use $\eta_{\frac{i}{2^{n_r}}}:= \mju$ in the
bisection formula \eqref{E:MidPtU} to get
\begin{equation}
\label{E:Diag}
F_{\frac{i}{2^{n_r}}}=\mju F_{\frac{i-1}{2^{n_r}}}+
\lp{1- \mju}F_{\frac{i+1}{2^{n_r}}}.
\end{equation}
We make the following observation.
\begin{observation}
\label{O:ccd}
Note that
\[ ||F_{\frac{i-1}{2^{n_r}}}+F_{\frac{i+1}{2^{n_r}}}-2F_{\frac{i}{2^{n_r}}}||_{L_1}
= || \mju G- (1-\mju) G||_{L_1} = \frac12 ||G||_{L_1} = 2^{-n_r},
\]
i.e. at $t = i/2^{n_r}$ there is a $c$-change of direction (cf \eqref{E:cChDir})
  in $F$ with $c= 1$.
  \end{observation}

\subsection{Details}

 First we use Lemma \ref{L:Split} to split $G_1:=F_1-F_0$ and
denote the obtained function by $\mju_{1}$ (it is clear that
$\mju_{1}(x)=x$, so this step is somewhat trivial, we discuss it
in detail for uniformity as it is quite similar to all further
steps). Denote by $D_R$ the centered at $0$ disc of radius $R$ in
$\mathbb{C}$. Since $\mju_1$ is a polynomial, there exists
$k_1\in\mathbb{N}$ such that
\begin{equation}
\label{E:firstU} |\mju_{1}^{(m)}(x)| \leq 2^{k_1-1}\omega_m\mbox{
and } |(1-\mju_{1})^{(m)}(x)| \leq 2^{k_1-1}\omega_m~~\forall
m\in\mathbb{N}~\forall x\in[0,1],
\end{equation}
and
\begin{equation}
\label{E:firstUC}\max\{|\mju_{1}(z)|,~ |(1-\mju_1)(z)|\} \leq
2^{k_1-1}~~\forall z\in D_1.
\end{equation}

We then use \eqref{E:MidPtDef} to define $F_t$ for the  first
$k_1$ generations of binary points. Note that as $F_0={\bf 0}$ and
$F_1={\bf 1}$, we get for all $1\leq n\leq k_1$, $1\leq i \leq
2^n$  that
\[ F_{\frac{i}{2^n}} (x) = \frac{i}{2^n} \ \ \mbox{ for all } \ x \in [0, 1], \]
 therefore
\begin{equation}
\label{E:Gen1}
 \ F_{\frac{i}{2^n}} (x) -F_{\frac{i-1}{2^n}} (x)
= \frac{1}{2^n} \ \mbox{ for all } \ 1\leq i \leq 2^n \hbox{ and
}x\in[0,1].
\end{equation}
For generation $n_1:=(k_1+1)$ we define $\disp
\eta_{\frac{i}{2^{n_1}}} := \mju_{1}$ and
\[ F_{\frac{i}{2^{n_1}}} = \eta_{\frac{i}{2^{n_1}}} F_{\frac{i-1}{2^{n_1}}}
+ \lp{1- \eta_{\frac{i}{2^{n_1}}}} F_{\frac{i+1}{2^{n_1}}} \ \
\mbox{ for every odd }  \ \ 1 \leq i \leq 2^{n_1}-1.\] Since
\[\ F_{\frac{i}{2^{n_1}}} - F_{\frac{i-1}{2^{n_1}}} = \lp{1- \eta_{\frac{i}{2^{n_1}}}}
\lp{F_{\frac{i+1}{2^{n_1}}} - F_{\frac{i-1}{2^{n_1}}}} \] and
\[  F_{\frac{i+1}{2^{n_1}}} - F_{\frac{i}{2^{n_1}}} = \eta_{\frac{i}{2^{n_1}}}
\lp{F_{\frac{i+1}{2^{n_1}}} - F_{\frac{i-1}{2^{n_1}}}}, \] we
obtain from \eqref{E:firstU} that for all $m\in\mathbb{N}$
\begin{equation}
\label{E:Diff1} |F^{(m)}_{\frac{i}{2^{n_1}}}(x) -
F^{(m)}_{\frac{i-1}{2^{n_1}}}(x)| \leq
2^{k_1-1}\omega_m\frac{1}{2^{k_1}} = \frac{\omega_m}{2} \ \ \mbox{
for every } \ \ 1\leq i \leq 2^{n_1} \hbox{ and }x\in[0,1],
\end{equation}
and from \eqref{E:firstUC} that
\begin{equation}
\label{E:C1} |F_{\frac{i}{2^{n_1}}}(z) -
F_{\frac{i-1}{2^{n_1}}}(z)| \leq 2^{k_1-1}\frac{1}{2^{k_1}} =
\frac{1}{2} \ \ \mbox{ for every } z\in D_1.
\end{equation}

In the next stage of the construction we use Lemma \ref{L:Split}
to split all of the $2^{k_1+1}=2^{n_1}$ functions  $\disp
G_{\frac{i}{2^{n_1}}}:= F_{\frac{i}{2^{n_1}}} -
F_{\frac{i-1}{2^{n_1}}} $ with $\disp 1 \leq i \leq 2^{n_1}$. We
denote the obtained $2^{n_1}$ functions by $\mju_{2,i}$, $1\le
i\le 2^{n_1}$.\medskip

Since $G_{\frac{i}{2^{n_1}}}$ and $\mju_{2,i}$ are polynomials, we
can pick $k_2\in\mathbb{N}$, such that for all
$i\in\{1,\dots,2^{n_1}\}$
\begin{equation}
\label{E:secondU}
\max\{|(\mju_{2,i}G_{\frac{i}{2^{n_1}}})^{(m)}(x)|,
|((1-\mju_{2,i})G_{\frac{i}{2^{n_1}}})^{(m)}(x)|\} \leq
2^{k_2-2}\omega_m~~\forall m\in\mathbb{N}~\forall x\in[0,1],
\end{equation}
and
\begin{equation}
\label{E:secondUC} \max\{|(\mju_{2,i}G_{\frac{i}{2^{n_1}}})(z)|,
|((1-\mju_{2,i})G_{\frac{i}{2^{n_1}}})(z)|\} \leq
2^{k_2-2}~~\forall z\in D_2.
\end{equation}

Then we use \eqref{E:MidPtDef} to define $F_t$ for the  next $k_2$
generations of binary points. Let $n_2:=n_1+k_2+1$. We obtain from
\eqref{E:secondU} that for all $m\in\mathbb{N}$
\begin{equation}
\label{E:Diff2} |F^{(m)}_{\frac{i}{2^{n_2}}}(x) -
F^{(m)}_{\frac{i-1}{2^{n_2}}}(x)| \leq
2^{k_2-2}\omega_m\frac{1}{2^{k_2}} = \frac{\omega_m}{2^2} \ \
\mbox{ for every } \ \ 1\leq i \leq 2^{n_2} \hbox{ and }x\in[0,1],
\end{equation}
and from \eqref{E:secondUC} that
\begin{equation}
\label{E:C2} |F_{\frac{i}{2^{n_2}}}(z) -
F_{\frac{i-1}{2^{n_2}}}(z)| \leq 2^{k_2-2}\frac{1}{2^{k_2}} =
\frac{1}{2^2} \ \ \mbox{ for every } z\in D_2.
\end{equation}

Assuming the construction was done up to ``$c$-changing
direction'' generation $n_{r-1}: = k_1+\dots + k_{r-1} + r-1$, we
define $k_r$ as follows. We use Lemma~\ref{L:Split} to split
$$\disp G_{\frac{i}{2^{n_{r-1}}}}:=
F_{\frac{i}{2^{n_{r-1}}}} - F_{\frac{i-1}{2^{n_{r-1}}}}, \ \ \
1\leq i \leq  2^{n_{r-1}},$$ denoting the obtained functions
$\mju_{r,i}$ and finding $k_{r} \geq 1$ such that for all $1\le
i\le 2^{n_{r-1}}$,
\begin{equation}
\label{E:rU}
\max\lc{\av{(\mju_{r,i}G_{\frac{i}{2^{n_{r-1}}}})^{(m)}(x)},
\av{\lp{(1-\mju_{r,i})G_{\frac{i}{2^{n_{r-1}}}}}^{(m)}(x)}} \leq
2^{k_{r}-r}\omega_m~~\forall m\in\mathbb{N}~\forall x\in[0,1],
\end{equation}
and
\begin{equation}
\label{E:rUC}
\max\lc{\av{(\mju_{r,i}G_{\frac{i}{2^{n_{r-1}}}})(z)},
\av{((1-\mju_{r,i})G_{\frac{i}{2^{n_{r-1}}}})(z)}} \leq
2^{k_r-r}~~\forall z\in D_r.
\end{equation}

Then we use \eqref{E:MidPtDef} to define $F_t$ for the  next $k_r$
generations of binary points. Let $n_r:=n_{r-1}+k_r+1$. We obtain
from \eqref{E:rU} that for all $m\in\mathbb{N}$
\begin{equation}
\label{E:Diffr} |F^{(m)}_{\frac{i}{2^{n_r}}}(x) -
F^{(m)}_{\frac{i-1}{2^{n_r}}}(x)| \leq
2^{k_r-r}\omega_m\frac{1}{2^{k_r}} = \frac{\omega_m}{2^r} \ \
\mbox{ for every } \ \ 1\leq i \leq 2^{n_r} \hbox{ and }x\in[0,1],
\end{equation}
and from \eqref{E:rUC} that
\begin{equation}
\label{E:Cr} |F_{\frac{i}{2^{n_r}}}(z) -
F_{\frac{i-1}{2^{n_r}}}(z)| \leq 2^{k_r-r}\frac{1}{2^{k_r}} =
\frac{1}{2^r} \ \ \mbox{ for every } z\in D_r.
\end{equation}

We continue in an obvious way and define $F_t$ for all binary
$t\in[0,1]$. As is mentioned, this defines $F_t\in L_1[0,1]$ for
all $t\in [0,1]$.

Our next goal is to show that the condition {\bf (A)} holds
following the plan described in Section \ref{S:Plan}.\medskip

So for each point $t\in[0,1]$ we need to find a sequence
$\{p(r,t)\}_{r=1}^\infty$ of binary points converging to $t$ and
satisfying the conditions in Section \ref{S:Plan}. We do this as
follows.

We let $p(1,t)$ to be the closest to $t$ binary number which
corresponds to a constant function (this can be non-unique only if
$t$ is binary of generation $n_1$, the first ``$c$-changing
direction'' generation). Here and later we resolve ties
arbitrarily.

We let $p(2,t)$ be the closest to $t$ binary number of generation
at most $n_2-1$ (this can be non-unique only if $t$ is binary of
generation $n_2$).

So on, we let $p(r,t)$ be the closest to $t$ binary number of
generation at most $n_r-1$. In this way we define $p(r,t)$ for all
$r\in\mathbb{N}$.
\medskip

It suffices to show that for all $r\geq 1$
\[|F_{p(r+1,t)}^{(m)}(x)-F_{p(r,t)}^{(m)}(x)|\le\frac{\omega_m}{2^r}\]
for all $m \in \mathbb{N}$ and all $x\in[0,1]$, and
\[|F_{p(r+1,t)}(z)-F_{p(r,t)}(z)|\le\frac{1}{2^r}\quad\forall z\in D_r.\]
These conditions follow from \eqref{E:Diffr} and \eqref{E:Cr},
respectively, since $F_{p(r+1,t)}-F_{p(r,t)}$ is a multiple of
$F_{\frac{i}{2^{n_{r}}}}-F_{\frac{i-1}{2^{n_{r}}}}$ with
coefficient in $[0,1]$. Therefore the condition {\bf (A)} is
satisfied.

\subsection{Nondifferentiability}

\begin{lemma}\label{L:cCD_implND} Suppose that a map $F:[0,1]\to L_1[0,1]$ is such that
there is $c>0$ and an infinite sequence $\{n_r\}_{r=1}^\infty$ of
positive integers such that the condition \eqref{E:cChDir} is
satisfied for each $n=n_r$ and each odd natural number $i\le 2^{n_r}$.
Then the map $F$ does not have points of differentiability.
\end{lemma}

This argument is essentially known, see \cite[Theorem 5.21, (iii)
$\Rightarrow$ (i)]{BL00}. For convenience of the reader we present
it.

\begin{proof}
Assume that $t_0\in[0,1]$ is a point of differentiability of $F_t$
(with respect to $t$) and let $D\in L_1[0,1]$ be the corresponding
derivative. Then
\begin{equation}
\label{E:diff}
 \nm{F_{t_0+h} -F_{t_0} -hD}_{L_1} = o(h) \ \ \mbox{ as $|h|\downarrow 0$ and  } t_0+h \in [0,1].
\end{equation}
Since for every $n$ fixed
\[ \bigcup_{i \mbox{ \footnotesize odd}} \ls{ \frac{i-1}{2^n}, \frac{i+1}{2^n}}=[0,1], \]
there exists a sequence of binary numbers in canonical form
$i_n/2^n$ such that
 \[ t_0\in \ls{ \frac{i_n-1}{2^n}, \frac{i_n+1}{2^n}}  \mbox{ for every } n.\]
Condition \eqref{E:diff} implies
\[ 2^{n}\nm{F_\frac{i_n\pm1}{2^{n}} -F_{t_0} - \lp{\frac{i_n\pm1}{2^{n}}-t_0}D}_{L_1} \to 0
\ \ \mbox{ as } \ \  n \to \infty,\]
and
\[ 2^{n}\nm{F_\frac{i_n}{2^{n}} -F_{t_0} - \lp{\frac{i_n}{2^{n}}-t_0}D}_{L_1} \to 0
\ \ \mbox{ as } \ \ n \to \infty.\]
By the triangle inequality we obtain
\[ 2^{n}\nm{F_\frac{i_n+1}{2^{n}}  +F_\frac{i_n - 1}{2^{n}} - 2F_\frac{i_n}{2^{n}} }_{L_1} \to 0
 \ \ \mbox{ as } \ \ n \to \infty.\]
However, this contradicts the fact that $F$ ``$c$-changes
direction'' at every $\disp \frac{i_{n_r}}{2^{n_r}}$.
\end{proof}

\section{Acknowledgement}

The authors would like to thank the referee for the
suggestion to state and to prove the main result in a stronger
form. The second named author gratefully acknowledges the support
by the National Science Foundation grants DMS-1201269 and
DMS-1700176, and by the Summer Support of Research program of St.
John's University during different stages of work on this paper.
\end{large}

\begin{small}

\renewcommand{\refname}{
\end{small}


\end{document}